\documentclass[12pt]{article}
\usepackage{amsfonts,amssymb,epsfig}
\usepackage{amsmath}

\date{}
\begin{document}
\newtheorem{df}{Definition}
\newtheorem{thm}{Theorem}
\newtheorem{lm}{Lemma}
\newtheorem{pr}{Proposition}
\newtheorem{co}{Corollary}
\newtheorem{re}{Remark}
\newtheorem{note}{Note}
\newtheorem{claim}{Claim}
\newtheorem{problem}{Problem}

\def\R{{\mathbb R}}

\def\E{\mathbb{E}}
\def\calF{{\cal F}}
\def\N{\mathbb{N}}
\def\calN{{\cal N}}
\def\calH{{\cal H}}
\def\n{\nu}
\def\a{\alpha}
\def\d{\delta}
\def\t{\theta}
\def\e{\varepsilon}
\def\t{\theta}
\def\pf{ \noindent {\bf Proof: \  }}
\def\trace{\rm trace}
\newcommand{\qed}{\hfill\vrule height6pt
width6pt depth0pt}
\def\endpf{\qed \medskip} \def\colon{{:}\;}
\setcounter{footnote}{0}

\def\Lip{{\rm Lip}}

\renewcommand{\qed}{\hfill\vrule height6pt  width6pt depth0pt}

\title{No greedy bases for matrix spaces with mixed $\ell_p$ and $\ell_q$ norms
}

\author{Gideon Schechtman\thanks{Supported in part by the Israel Science Foundation. } } \maketitle

\begin{abstract}
We show that non of the spaces $(\bigoplus_{n=1}^\infty\ell_p)_{\ell_q}$, $1\le p\not= q<\infty$  have a greedy basis. This solves a problem raised by Dilworth, Freeman, Odell and Schlumprect.
Similarly, the spaces $(\bigoplus_{n=1}^\infty\ell_p)_{c_0}$, $1\le p<\infty$, and $(\bigoplus_{n=1}^\infty c_o)_{\ell_q}$, $1\le q<\infty$, do not have greedy bases. It follows from that and known results that a class of  Besov spaces on $\R^n$ lack greedy bases as well.
\end{abstract}

{\section{Introduction}

Given a (say, real) Banach space $X$ with a Schauder basis $\{x_i\}$, an $x\in X$ and an $n\in\mathbb{N}$ it is useful to determine the best $n$-term approximation to $x$ with respect to the given basis. I.e., to find a set $A\subset\mathbb{N}$ with $n$ elements and coefficients $\{a_i\}_{i\in A}$ such that
\[ 
\|x-\sum_{i\in A}a_ix_i\|=\inf\{\|x-\sum_{i\in B}b_ix_i\|; |B|=n, b_i\in \mathbb{R}\}
\]
or, given a $C<\infty$, at least to find such an $A\subset\mathbb{N}$ and coefficients $\{a_i\}_{i\in A}$ with
\[ 
\|x-\sum_{i\in A}a_ix_i\|\le C\inf\{\|x-\sum_{i\in B}b_ix_i\|; |B|=n, b_i\in \mathbb{R}\}.
\]
This problem attracted quite an attention in modern Approximation Theory. Of course one would also like to have a simple algorithm to find such a set $\{a_i\}_{i\in A}$. It would be nice if we could take $\{a_i\}_{i\in A}$ to be just the set of the $n$ largest, in absolute value, coefficients in the expansion of $x$ with respect to the basis $\{x_i\}$. Or, if this set is not unique, any such set. The basis $\{x_i\}$ is called {\em Greedy} if for some $C$ this procedure works; i.e., for all $x=\sum_{i=1}^\infty a_ix_i$, all $n\in\mathbb{N}$ and all $A\subset \mathbb{N}$, $|A|=n$, satisfying $\min\{|a_i|; i\in A\}\ge \max\{|a_i|; i\notin A\}$,
\[ 
\|x-\sum_{i\in A}a_ix_i\|\le C\inf\{\|x-\sum_{i\in B}b_ix_i\|; |B|=n, b_i\in \mathbb{R}\}.
\]

Konyagin and Temlyanov \cite{kt} provided a simple criterion to determine whether a basis is greedy: $\{x_i\}$ is greedy if and only if it is {\em unconditional} and {\em democratic}.

Recall that $\{x_i\}$ is said to be unconditional provided, for some $C<\infty$, all eventually zero coefficients $\{a_i\}$ and all sequences of signs $\{\e_i\}$,
\[
\|\sum\e_ia_ix_i\|\le C \|\sum a_ix_i\|.
\]
$\{x_i\}$ is said to be democratic provided for some $C<\infty$ and all finite $A,B\subset\mathbb{N}$ with $|A|=|B|$,
\[
\|\sum_{i\in A}x_i\|\le C\|\sum_{i\in B}x_i\|.
\]

We refer to \cite{dfos} for a survey of what is known about space that have or do not have greedy bases. In \cite{dfos} Dilworth, Freeman, Odell and Schlupmrecht determined which of the spaces $X=(\bigoplus_{n=1}^\infty\ell_p^n)_{\ell_q}$, $1\le p\not= q\le\infty$ (with $c_0$ replacing $\ell_\infty$ in case $q=\infty$) have a greedy basis. It turns out that this happens exactly when $X$ is reflexive. They also raise the question  of whether $(\bigoplus_{n=1}^\infty\ell_p)_{\ell_q}$, $1< p\not= q<\infty$ have greedy bases. Here we show that these spaces (as well as their non-reflexive counterparts) do not have greedy bases. By the Konyagin-Temlyanov characterization it is enough to prove that each unconditional basis of $(\bigoplus_{n=1}^\infty\ell_p)_{\ell_q}$, $1\le p\not= q\le\infty$ (with $c_0$ replacing $\ell_\infty$ in case $p$ or $q$ are $\infty$) has two subsequences, one equivalent to the unit vector basis of $\ell_p$ ($c_0$ if $p=\infty$) and one to the unit vector basis of $\ell_q$ ($c_0$ if $q=\infty$).

\begin{thm}\label{thm:main}
Each normalized unconditional basis of the spaces $(\bigoplus_{n=1}^\infty\ell_p)_{\ell_q}$, $1\le p\not= q<\infty$ has a subsequence equivalent to the unit vector basis of $\ell_p$ and another one equivalent to the unit vector basis of $\ell_q$.
Similarly, each unconditional basis of the spaces $(\bigoplus_{n=1}^\infty\ell_p)_{c_0}$, $1\le p<\infty$ (resp. $(\bigoplus_{n=1}^\infty c_o)_{\ell_q}$, $1\le q<\infty$) has a subsequence equivalent to the unit vector basis of $\ell_p$ (resp. $c_0$) and another one equivalent to the unit vector basis of $c_o$ (resp. $\ell_q$).
Consequently, none of these spaces have a greedy basis.
\end{thm}

For $1\le p,q<\infty$ the spaces $(\bigoplus_{n=1}^\infty\ell_p)_{\ell_q}$ are isomorphic to certain Besov spaces on $\R^n$. We refer to \cite{me} for the definition of the Besov spaces $B_p^{s,q}$ and for the fact that they are isomorphic to $(\bigoplus_{n=1}^\infty\ell_p)_{\ell_q}$. See in particular \cite[Section 6.10, Proposition 7]{me} (and also \cite[Section 2.9, Proposition 4]{me}). We thank P. Wojtaszczyk for this reference.

\begin{co}\label{co:besov}
Let $1\le p\not=q<\infty$ and $s$ any real number then the space $B_p^{s,q}$ does not have a greedy basis.
\end{co}
Recall that this stand in contrast with the main result in \cite{dfos} which states that, in the reflexive cases, the corresponding Besov spaces on $[0,1]$ do have greedy bases.

In the special case of $1<q<\infty$ and $p=2$ the theorem above was actually proved in \cite{sc}. There the isomorphic classification of the span of unconditional basic sequences in $(\bigoplus_{n=1}^\infty\ell_2)_{\ell_q}$, $1<q<\infty$, which span complemented subspaces were characterize. Although it is not stated there, the proof actually established the theorem above in these special cases. Shortly after \cite{sc} appeared Odell \cite{sc} strengthened the result and classified {\em all} the complemented subspace of $(\bigoplus_{n=1}^\infty\ell_2)_{\ell_q}$ (thus there is no wonder that \cite{sc} was forgotten...). We remark in passing that this special case of $p=2$ was of particular interest since $(\bigoplus_{n=1}^\infty\ell_2)_{\ell_q}$ is isomorphic to a complemented subspace of $L_q[0,1]$.

The first step in the proof in \cite{sc} is to reduce the case of a general unconditional basic sequence in $(\bigoplus_{n=1}^\infty\ell_2)_{\ell_q}$ whose span is complemented to one which is also a block basis of the natural basis of $(\bigoplus_{n=1}^\infty\ell_2)_{\ell_q}$. This reduction no longer hold for $p\not= 2$. The complications in the present note stems from this fact. The way we overcome it is by transferring the problem to a larger space (of arrays $\{a_{i,j,k}\}$) of mixed $q,p$ and 2 norms. Unfortunately, this makes the notations quite cumbersome.

\section{Preliminaries}

$Z_{q,p}$, $1\le p,q<\infty$ will denote here the space of all matrices $a=\{a(i,j)\}_{i,j=1}^\infty$ with norm
\[
\|a\|=\|a\|_{q,p}=(\sum_{j=1}^\infty(\sum_{i=1}^\infty|a(i,j)|^p)^{q/p})^{1/q}.
\]
If $p$ or $q$ are $\infty$ we replace the corresponding $\ell_p$ or $\ell_q$ norm by the $\ell_\infty$ norm and continue to denote by $Z_{q,p}$ the completion of the space of finitely supported matrices under this norm. (Thus, $c_0$ replacing $\infty$ would be a more precise notation in this case but, since it would complicated our statements, we prefer the above notation.)
The spaces $Z_{q,p}$ are the subject of investigation of this paper. They are more commonly denoted by $\ell_q(\ell_p)$ or $(\bigoplus_{n=1}^\infty \ell_p)_{\ell_q}$ (as we have done in the introduction) but since we shall be forced to also consider more complicated spaces with mixed norms we prefer the notation above.

If $\{k_n\}_{n=1}^\infty$ is any sequence of positive integers, we shall denote by
$Z_{q,p;\{k_n\}}$, the subspace of $Z_{q,p}$ consisting of matrices $a$ satisfying $a(i,j)=0$ for all $i>k_j$.

We also denote by $Z_{q,p,r}$ (we'll use this only for $r=2$) the spaces of arrays $a=\{a(u,i,j)\}_{u,i,j=1}^\infty$ with norm
\[
\|a\|=\|a\|_{q,p,r}=(\sum_{j=1}^\infty(\sum_{i=1}^\infty(\sum_{u=1}^\infty|a(u,i,j)|^r)^{p/r})^{q/p})^{1/q},
\]
with the same convention as above when one of $p,q$ (or $r$) is $\infty$.
Similarly, $Z_{q,p;\{k_n\},r}$ denotes the subspace of $Z_{q,p,r}$ consisting of arrays $a$ satisfying $a(u,i,j)=0$ for all $i>k_j$.

By $P_n$ we denote the natural projection onto the $n$-th column in $Z_{q,p}$; i.e, $P_n(\{a(i,j)\})=\{\bar a(i,j)\}$, where $\bar a(i,j)=a(i,j)$ if $j=n$  and $\bar a(i,j)=0$ otherwise. Similarly, $P_n^k$ denotes the natural projection onto the first $k$ elements in the $n$-th column.
$Q_N$ denotes $\sum_{n=1}^NP_n$.

Given a Banach lattice $X$, an $1< r<\infty$ and $x_1,x_2,\dots\in X$ one can define the operation $(\sum|x_n|^r)^{1/r}$ in a manner consistent with what we usually mean by such an operation (when $X$ is a lattice of functions or sequences, for example). See e.g. \cite[Section 1.d]{lt2} for this and what follows.

In particular if $X$ has a 1-unconditional basis $\{e_i\}$ (which is the only kind of lattices we'll consider here) then for  $x_n=\sum_{i=1}^\infty a_i^ne_i$, $n=1,2,\dots,N$, $(\sum|x_n|^r)^{1/r}=\sum_{i=1}^\infty(\sum_{n=1}^N|a_i^n|^r)^{1/r}e_i$.

Recall that $X$ is said to be $r$-convex (resp. $r$-concave) with constant $K$ if for all $n$ and all $x_1,x_2,\dots,x_n\in X$
\[
\|(\sum_{i=1}^n|x_i|^r)^{1/r}\|\le K (\sum_{i=1}^n\|x_i\|^r)^{1/r} \
(\hbox{resp.} \ (\sum_{i=1}^n\|x_i\|^r)^{1/r}\le K\|(\sum_{i=1}^n|x_i|^r)^{1/r}\|).
\]
$X$ is said to be $r$-convex (resp. $r$-concave) if it is $r$-convex (resp. $r$-concave) with some constant $K<\infty$. $Z_{q,p}$ is easily seen to be $\min\{p,q\}$-convex with constant $1$ and $\max\{p,q\}$-concave with constant $1$.

It is also known that $X$ is $r$-convex (resp. $r$-concave) if and only if its dual $X^*$ is $r^\prime$-concave (resp. $r^\prime$-convex) where $r^\prime=r/(r-1)$.

Given a Banach lattice $X$ we denote by $X(\ell_2)$ the (completion of the) space of (finite) sequences $x=(x_1,x_2,\dots)$ of elements of $X$ for which the norm
\[
\|x\|=\|(\sum|x_j|^2)^{1/2}\|
\]
is finite. If $X$ has a $1$-unconditional basis $\{e_j\}$ then this is just the (completion of the) space of matrices $a=\{a(i,j)\}$ (with only finitely many non-zero entries) with norm
\[
\|a\|=\|\sum_{j=1}^\infty(\sum_{i=1}^\infty |a(i,j)|^2)^{1/2}e_i\|.
\]

The following two lemmas are well known but maybe hard to find so we reproduce their proofs.

\begin{lm}\label{lm:l_p}
Let $\{x_i\}_{i=1}^\infty$ be a normalized unconditional basic sequence in $Z_{q,p}$, $1\le p<q\le\infty$. If for some $\e>0$ and $N\in \mathbb{N}$  $\|Q_Nx_i\|>\e$ for all $i$ then $\{x_i\}_{i=1}^\infty$ has a subsequence equivalent to the unit vector basis of $\ell_p$.
\end{lm}

\pf
Assume first $p>1$. Given a sequence of positive $\e_i$-s and passing to a subsequence (which without loss of generality we assume is the all sequence) we can assume that there is a sequence of $\{y_i\}$ of vectors disjointly supported with respect to the natural basis of $Z_{q,p}$ such that $\|x_i-y_i\|<\e_i$ for all $i$. (Use the fact that $\{x_i\}$ doesn't have a subsequence equivalent to the unit vector basis of $\ell_1$ and the argument for Proposition 1.a.12 in \cite{lt1}, for example). $\{y_i\}$ is 1-dominated by the unit vector basis of $\ell_p$ and dominates $\{Q_Ny_i\}$ which in turn $C$-dominates the unit vector basis of $\ell_p$ for $C=1/(\e-\sup\e_i)$; i.e,
\[
(\sum_{i=1}^\infty |a_i|^p)^{1/p}\ge \|\sum_{i=1}^\infty a_iy_i\|\ge (\e-\sup\e_i)(\sum_{i=1}^\infty |a_i|^p)^{1/p}
\]
for all scalars $\{a_i\}$. If the $\e_i$-s are small enough a similar inequality holds for the (sub)sequence $\{x_i\}$.

If $p=1$ then given a sequence of positive $\e_i$-s and passing to a subsequence (which without loss of generality we assume is the all sequence) we can assume that there is a vector $y$ and sequence of $\{y_i\}$ of vectors all disjointly supported with respect to the natural basis of $Z_{q,p}$ such that $\|x_i-y-y_i\|<\e_i$ for all $i$. If $y\not=0$ and the $\{\e_i\}$ are small enough then, using the unconditionality} $\{x_i\}$ is clearly equivalent to the unit vector basis of $\ell_1$. If $y=0$ the same argument as for $p>1$ works here too.
\endpf

\begin{lm}\label{lm:complementation}
Let $\{x_i\}$ be a $K$-unconditional basic sequence in a Banach lattice which is $r$-concave for some $r<\infty$
Let $\bar x_i\in X(\ell_2)$ be defined by $(0,\dots,0,x_i,0,\dots)$, $x_i$ in the $i$-th place. Then the sequences $\{x_i\}$ in $X$ and $\{\bar x_i\}$ in $X(\ell_2)$ are equivalent.

If in addition $X$ is also $s$-convex for some $s>1$ and $[x_i]$, the closed linear span of $\{x_i\}$, is complemented in $X$ then $[\bar x_i]$ is complemented in $X(\ell_2)$.
\end{lm}

\pf
The first assertion, due to Maurey, can be found in \cite{ma} or \cite[Theorem1.d.6(i)]{lt2}.
The second is probably harder to find so we reproduce it.
Let $P=\sum_{i=1}^\infty x_i^*\otimes x_i$, with $x^*_i\in X^*$, be the projection onto $[x_i]$; i.e.,
\[
P(x)=\sum_{i=1}^\infty x_i^*(x) x_i\ \ x\in X.
\]
Define $\bar P=\sum_{i=1}^\infty \bar x_i^*\otimes \bar x_i$ ($\bar x^*_i\in X^*(\ell_2)=X(\ell_2)^*$); i.e.,
\[
\bar P(x)=\sum_{i=1}^\infty \bar x_i^*(x) \bar x_i\ \ x\in  X(\ell_2).
\]
Using the facts that $\{\bar x_i\}$ is equivalent to $\{ x_i\}$, $\{\bar x_i^*\}$ is equivalent to $\{ x_i^*\}$, and $\{\bar x_i^*,\bar x_i\}$ is a biorthogonal sequence, it is easy to see that $\bar P$ is a bounded projection on $X(\ell_2)$ with range $[\bar x_i]$.
\endpf

\section{Proof of the main result, the reflexive case}

Since the non-reflexive cases (i.e., when $p$ or $q$ are $1$ or $\infty$) of Theorem \ref{thm:main} require a bit different treatment and since the problem raised in \cite{dfos} was restricted to the reflexive cases only, we prefer to delay the proof of the non-reflexive cases to the next section.

\begin{pr}\label{pr:main1}
Let $\{x_i\}_{i=1}^\infty$ be a normalized unconditional basic sequence in $Z_{q,p}$, $1<p,q<\infty$ such that $[x_i]_{i=1}^\infty$ is complemented in $Z_{q,p}$. If no subsequence of $\{x_i\}_{i=1}^\infty$ is equivalent to the unit vector basis of $\ell_p$ then $[x_i]_{i=1}^\infty$ isomorphically embeds in $Z_{q,p;\{n\},2}$ as a complemented subspace.
\end{pr}

\pf We may clearly assume $p\not= q$ and by duality that $q>p$. Let $\{\e_n\}_{n=1}^\infty$ be a sequence of positive numbers. By Lemma \ref{lm:l_p} for all $n$ only finitely many of the $x_i$-s satisfy $\|P_nx_i\|\ge\e_n$. Consequently, for each $n\in\mathbb{N}$ there is a $k_n\in\mathbb{N}$ such that $\|(P_n-P_n^{k_n})x_i\|<\e_n$ for all $i$. We denote
$Q=\sum_{n=1}^\infty P_n^{k_n}$.
In the case $p=2$ we showed in \cite{sc} that without loosing generality we can assume that $\{x_i\}$ is a block basis of the natural basis of $Z_{q,p}$ and then $\{Qx_i\}$ and $\{(I-Q)x_i\}$ are also unconditional basic sequences. This is no longer true when $p\not= 2$. We overcome this difficulty by switching to the larger space $Z_{q,p,2}$. Define for each $i$  $\bar x_i\in Z_{q,p,2}$ by
\[
\bar x_i(w,u,v)=\left\{
                   \begin{array}{ll}
                    x_i(u,v), & \hbox{if} \ w=i; \\
                     0, & \hbox{if} \ w\not= i.
                   \end{array}
                 \right.
\]
Let the projection $P$ from $Z_{q,p}$ onto $[x_i]$ be given by
\[
Px=\sum_{i=1}^\infty x^*_i(x)x_i
\]
where $\{x^*_i\}$ in $Z_{q^\prime,p^\prime}=Z^*_{q,p}$ ($p^\prime=p/(p-1), \ q^\prime=q/(q-1)$) satisfy $x^*_i(x_j)=\delta_{i,j}$, $i,j=1,2,\dots$. Then by Lemma \ref{lm:complementation}
\[
\bar Px=\sum_{i=1}^\infty \bar x^*_i(x)\bar x_i
\]
is a bounded projection from $Z_{q,p,2}$ onto $[\bar x_i]$ and $\{x_i\}_{i=1}^\infty$ is equivalent to $\{\bar x_i\}_{i=1}^\infty$.

We denote by $\bar P_n=P_n\otimes I_{\ell_2}$ on $Z_{q,p,2}$; i.e, $\bar P_n(x)(w,u,v)=P_n(x(w,\cdot,\cdot))(u,v)$. We also similarly denote $\bar P_n^k=P_n^k\otimes I_{\ell_2}$, $\bar Q_N=Q_N\otimes I_{\ell_2}$, and $\bar Q=Q\otimes I_{\ell_2}$. Note that now $\{\bar Q\bar x_i\}$ and $\{(I-\bar Q)\bar x_i\}$ are also unconditional basic sequences.
We would like to show that if $\e_n\to 0$ fast enough, then $\{\bar Q\bar x_i\}$ is equivalent to $\{\bar x_i\}$ and thus to $\{x_i\}$ and that $[\bar Q\bar x_i]$ is complemented.

Now,
\[
(I-\bar Q)\bar P\bar Q\bar x_n=\sum_{i=1}^\infty \bar x_i^*(\bar Q\bar x_n)(I-
\bar Q)\bar x_i, \ \ n=1,2,\dots.
\]
The operator $(I-\bar Q)\bar P$ sends the span of the unconditional basic sequence $\{\bar Q\bar x_n\}$ to the span of 
the unconditional basic sequence $\{(I-\bar Q)\bar x_n\}$ thus the diagonal operator $D$ defined by
\[
D\bar Q\bar x_n=\bar x_n^*(\bar Q\bar x_n)(I-
\bar Q)\bar x_n, \ \ n=1,2,\dots.\
\]
is bounded (see e.g. \cite{to} or
\cite[Proposition 1.c.8]{lt1}). If we show that $\bar x_n^*(\bar Q\bar x_n)$ are uniformly bounded away from zero this will show that $\{\bar Q\bar x_n\}$ dominates $\{(I-\bar Q)\bar x_n\}$ and thus also $\{\bar x_n\}=\{(I-\bar Q)\bar x_n+\bar Q\bar x_n\}$. That  $\{\bar Q\bar x_n\}$ is dominated by $\{\bar x_n\}$ is clear from the boundedness of $\bar Q$. This will show that $\{\bar Q\bar x_n\}$ is equivalent to $\{x_n\}$. To show that $\bar x_n^*(\bar Q\bar x_n)$ are uniformly bounded away from zero note that
\[
\bar x_n^*(\bar Q\bar x_n)= 1-\bar x_n^*((I-\bar Q) \bar x_n)
\]
and that
\[
|\bar x_n^*((I-\bar Q) \bar x_n)|\le \|\bar P(I-\bar Q) \bar x_n\|\le \|\bar P\|\sum_{i=1}^\infty \e_i.
\]
So, if $\|\bar P\|\sum_{i=1}^\infty \e_i<1/2$, then $\bar x_n^*(\bar Q\bar x_n)\ge 1/2$ for all $n$.

We still need to show that $[\bar Q\bar x_n]$ is complemented. Note that $\{\frac{\bar x_n^*}{\bar x_n^*(\bar Q\bar x_n)}, \bar Q\bar x_n\}$ is a biorthogonal sequence such that $\{\bar Q\bar x_n\}$ is equivalent to $\{\bar x_n\}$ and $\{\frac{\bar x_n^*}{\bar x_n^*(\bar Q\bar x_n)}\}$ is dominated by $\{x_n^*\}$. It follows that
\[
x\to\sum_{n=1}^\infty \frac{\bar x_n^*(x)}{\bar x_n^*(\bar Q\bar x_n)} \bar Q\bar x_n
\]
defines a bounded projection with range $[\bar Q\bar x_n]$.

We have shown that $[x_i]$ embeds complementably into $Z_{q,p;\{k_n\},2}$ for some sequence of positive integers $\{k_n\}$. This last space is clearly isometric to a norm one complemented subspace of $Z_{q,p;\{n\},2}$.
\endpf

In the case $p=2$ the argument above simplifies and actually shows that under the assumptions of Proposition \ref{pr:main1} we can strengthen the conclusion to: $[x_i]$ embeds complementably in $Z_{q,2;\{n\}}$ (which is isomorphic to $\ell_q$). We will not dwell on it farther as this is contained in \cite{sc}. The next proposition combained with the previous one will show in particular that any unconditional basis of $Z_{q,p}$ contains a subsequence equivalent to the unit vector basis of $\ell_p$. We'll need to use this also in the next section so we include the non-reflexive cases here as well.

\begin{pr}\label{pr:main2} Let $1\le p,q\le \infty$.
If $p\not=2,q$, then $\ell_p$ ($c_0$ in case $p=\infty$) does not embed into $Z_{q,p;\{n\},2}$.
\end{pr}

\pf
Assume $\ell_p$ or $c_0$ embeds into $Z_{q,p;\{n\},2}$. Passing to a subsequence of the image of the unit vector basis of $\ell_p$ or $c_0$, taking successive differences (this is needed only in the case $p=1$) and using a simple perturbation argument, we may assume that some normalized block basis $\{x_i\}$ of the natural basis of $Z_{q,p;\{n\},2}$ is equivalent to the unit vector basis of $\ell_p$ ($c_0$ if $p=\infty$). Let $P_{n,m}$, $m=1,2\dots$, $1\le n\le m$, denote the canonical projection onto the $n,m$ copy of $\ell_2$ in $Z_{q,p;\{n\},2}$:
\[
P_{n,m}(\{a(w,u,v)\})=\{\bar a(w,u,v)\}
\]
where
\[
\bar a(w,u,v)=\left\{
                   \begin{array}{ll}
                    a(w,u,v), & \hbox{if} \ u=n,\ v=m; \\
                     0, & \hbox{otherwise}.
                   \end{array}
                 \right.
\]
Assume first $p>2$. For each $n,m$ $P_{n,m}$ acts as a compact operator from $[x_i]$ to $\ell_2$ as every bounded operator from $\ell_p$, $p>2$ or $c_0$ to $\ell_2$ do. Consequently, given a sequence of positive numbers $\{\e_{n,m}\}$, we can find $k_{n,m}\in\mathbb{N}$ such that $\|(P_{n,m}-P_{n,m}^{k_{n,m}})_{|[x_i]}\|<\e_{n,m}$ for all $n,m$. 
Then, if $\sum_{n,m}\e_{n,m}$ is small enough 
\[
(\sum_{n,m}P_{n,m}^{k_{n,m}})_{|[x_i]}
\] 
is an isomorphism
and we get that $[x_i]$ embeds into $Z_{q,p;\{n\},2;\{k_{n,m}\}}$. Now for each finite $m$ and $k$ the $\ell_p^m$ sum of $\ell_2^k$-s 2-embeds into $\ell_p^N$ for some $N$ depending only on $p,m$ and $k$. It thus follows that $[x_i]$ embeds into $Z_{q,p;\{k_n\}}$ for some sequence of positive integers $\{k_n\}$. Passing to a farther subsequence of $\{x_i\}$, we get that the unit vector basis of $\ell_p$ (or $c_0$ in the case $p=\infty$) is equivalent to that of $\ell_q$ which is a contradiction.

The case $1\le p<2$ is just a bit more complicated. Here $P_{n,m}$ doesn't act as a compact operator from $[x_i]$ to $\ell_2$ but it is still strictly singular. Consequently, for each $n,m$ and $l$ we can find a normalised block basis of $\{x_i\}_{i=l}^\infty$ such that $\|(P_{n,m})_{|[x_i]_{i=l}^\infty}\|<\e_{n,m}$ and consequently there is a block basis of $\{x_i\}$ {\em whose first $l-1$ terms are just $x_1,\dots,x_{l-1}$}, and $k_{n,m,l}$ such that
\[
\|(P_{n,m}-P_{n,m}^{k_{n,m,l}})_{|[x_i]}\|<\e_{n,m}.
\]
A simple diagonal argument will now produce a normalised block basis $\{z_i\}$ of $\{x_i\}$ and natural numbers $k_{n,m}$-s such that
\[
(\sum_{n,m}P_{n,m}^{k_{n,m}})_{|[z_i]}
\]
an isomorphism. Since $\{z_i\}$ is equivalent to the unit vector basis of $\ell_p$ we get that $\ell_p$ embeds into $Z_{q,p;\{n\},2;\{k_{n,m}\}}$. The rest of the proof in this case is the same as in the case $p>2$.
\endpf

We are now aiming at proving that every normalized unconditional basis of $Z_{q,p}$ contains a subsequence equivalent to the unit vector basis of $\ell_q$.

\begin{pr}\label{pr:main3}Let $\{x_i\}_{i=1}^\infty$ be a normalized unconditional basic sequence in $Z_{q,p}$, $1<p,q<\infty$ such that $[x_i]_{i=1}^\infty$ is complemented in $Z_{q,p}$. If no subsequence of $\{x_i\}_{i=1}^\infty$ is equivalent to the unit vector basis of $\ell_q$ then $[x_i]_{i=1}^\infty$ isomorphically embeds in $Z_{p,2}$ as a complemented subspace.
\end{pr}

\pf We may assume $q<p$. We first claim that for each $\e>0$ there is an $N$ such that $\|(I-Q_N)x_i\|<\e$ for each $i=1,2,\dots$. Indeed if this is not the case then there is an $\e>0$, a sequence $0=N_1<N_2<\cdots$ in $\mathbb{N}$ and a subsequence $\{y_i\}$ of $\{x_i\}$ such that $\|(Q_{i+1}-Q_i)y_i\|\ge\e$ for all $i$. Passing to a further subsequence and a small perturbation we may assume that $\{y_i\}$ is a block basis of the natural basis of $Z_{q,p}$.
Then, since $q<p$, for all scalars $\{a_i\}$,
\[
(\sum_{i=1}^\infty |a_i|^q)^{1/q}\ge \|\sum_{i=1}^\infty a_iy_i\|\ge\|\sum_{i=1}^\infty a_i(Q_{i+1}-Q_i)y_i\|\ge
\e(\sum_{i=1}^\infty |a_i|^q)^{1/q}
\]
in contradiction to the fact that no subsequence of the $\{x_i\}$ is equivalent to the  unit vector basis of $\ell_q$.

The rest of the proof is now similar to that of Proposition \ref{pr:main1}, only a bit simpler. Fix an $\e>0$ and let $N$ be as in the beginning of this proof. Let $P=\sum_{i=1}^\infty x_i^*\otimes x_i$ be the projection onto $[x_i]$ and let $\{\bar x_i\}$ (in $Z_{q,p,2}$), $\bar P$ and $\bar Q_N$ be as in the proof of Proposition \ref{pr:main1}.  Consider the operator $(I-\bar Q_N)\bar P$ as acting from the span of the unconditional basic sequence $\{\bar Q_N \bar x_i\}$ to the span of the unconditional sequence $\{(I-\bar Q_N) \bar x_i\}$:
\[
(I-\bar Q_N)\bar P \bar Q_N\bar x_n=\sum_{i=1}^\infty \bar x_i^*(Q_N\bar x_n)(I-\bar Q_N)\bar x_i, \ n=1,2,\dots.
\]
Its diagonal defined by
\[
D\bar Q_N\bar x_n=\bar x_n^*(Q_N\bar x_n)(I-\bar Q_N)\bar x_n, \ n=1,2,\dots
\]
is bounded (\cite{to} or \cite{lt1}). So if we show that $\bar x_n^*(\bar Q_N\bar x_n)$ are bounded away from zero then the sequence $\{\bar Q_N \bar x_i\}$ will dominate the sequence $\{(I-\bar Q_N) \bar x_i\}$ and thus also $\{\bar x_i\}$ and $\{ x_i\}$. This will also show that
\[
x\to\sum_{n=1}^\infty \frac{\bar x_n^*(x)}{\bar x_n^*(\bar Q_N\bar x_n)} \bar Q_N\bar x_n
\]
defines a bounded projection from $\bar Q_NZ_{q,p,2}$ (which is isomorphic to $Z_{p,2}$) onto $[\bar Q_N\bar x_i]$ (which is isomorphic to $[ x_i]$).

To show that $\bar x_n^*(\bar Q_N\bar x_n)$ are bounded away from zero note that
\[
\bar x_n^*(\bar Q_N\bar x_n)= 1-\bar x_n^*((I-\bar Q_N) \bar x_n)
\]
and that
\[
|\bar x_n^*((I-\bar Q_N) \bar x_n)|\le \|\bar P(I-\bar Q_N) \bar x_n\|\le \|\bar P\|\e.
\]
So, if $\|\bar P\|\e<1/2$, then $\bar x_n^*(\bar Q\bar x_n)\ge 1/2$ for all $n$.
\endpf

\begin{re}
With a bit more effort one can strengthn the conclusion of Proposition \ref{pr:main3} to: $[x_i]_{i=1}^\infty$ is isomorphic to $\ell_p$. This is done by first showing that one can embed $[x_i]_{i=1}^\infty$ as a complemented subspace in $Z_{p,2;\{n\}}$ which is isomorphic to $\ell_p$ and using the fact that any infinite dimensional complemented subspace of $\ell_p$ is isomorphic to $\ell_p$. 
\end{re}

\medskip

\noindent{\bf Proof of Theorem \ref{thm:main} in the reflexive case:} Propositions \ref{pr:main1} and \ref{pr:main2} show that any normalized unconditional basis of $Z_{q,p}$, $1<p,q<\infty$, has a subsequence equivalent to the unit vector basis of $\ell_p$. To show that any such basis also has a subsequence equivalent to the unit vector basis of $\ell_q$ we need, in view of Proposition \ref{pr:main3}, only prove that $Z_{q,p}$ doesn't embed complementably into $Z_{p,2}$ for $1<q\not=p<\infty$.  This can probably be done directly (especially in the case $q\not=2$ in which case it is also true that $\ell_q$ does not embed into $Z_{p,2}$) but it also follows from the main theorems of \cite{sc} and \cite{od} in which the complemented subspaces of $Z_{p,2}$ (in \cite{sc} only those with unconditional basis) where characterized.
\endpf

\section{Proof of the main result, the non-reflexive case}

Recall that the subscript $\infty$ in $Z_{\infty,p}$ refers, by our convention, to the $c_0$ (rather than $\ell_\infty$) sum. Similatly, the subscript $\infty$ in $Z_{q,\infty}$ refers to the $q$ sum of $c_0$. We are going to show that any unconditional basis of each of the spaces $Z_{q,p}$, $p\not= q$, when at least one of $p$ or $q$ is $1$ or $\infty$ contains a subsequence equivalent to the unit vector basis of $\ell_p$ ($c_0$ if $p=\infty$) and another subsequence  equivalent to the unit vector basis of $\ell_q$ ($c_0$ if $q=\infty$). 

The spaces $Z_{1,\infty}$ and $Z_{\infty,1}$ (as well as $Z_{1,2}$ and $Z_{\infty,2}$) have unique, up to permutation, unconditional bases \cite{bclt}. These bases clearly contain a subsequence equivalent to the unit vector basis of $c_0$ and another one equivalent to the unit vector basis of $\ell_1$, so we only need to deal with the spaces $Z_{\infty,p}$, $1<p<\infty$, and their duals $Z_{1,p^\prime}$ and with $Z_{q,\infty}$, $1<q<\infty$, and their duals $Z_{q^\prime,1}$.

We shall need some classical results concerning unconditional bases and duality. These can be found conveniently in sections 1.b. and 1.c. of \cite{lt1}. $\ell_1$ does not isomorphically embed into  $Z_{\infty,p}$, $1<p<\infty$, (resp. into $Z_{q,\infty}$, $1<q<\infty$) (this follows for example from the fact that these spaces are $p$ (resp. $q$) convex). It thus follows from a theorem of James \cite{ja} or see \cite[Theorem 1.c.9]{lt1} that any unconditional basis of these spaces is shrinking. See \cite[Proposition 1.b.1]{lt1} for the the definition of a shrinking basis as well as for the fact that then the biorthogonal basis is an unconditional basis of the dual space $Z_{1,p^\prime}$, $1<p<\infty$, (resp. $Z_{q^\prime,1}$, $1<q<\infty$). Thus, in order to prove Theorem \ref{thm:main} in the non-reflexive cases, if would be enough to show that any normalized unconditional basis of  $Z_{1,p}$, $1<p<\infty$, (resp. $Z_{q,1}$, $1<q<\infty$) contains a subsequence equivalent to the unit vector basis of $\ell_1$ and another subsequence equivalent to the unit vector basis of $\ell_p$ (resp. $\ell_q$).

Let $\{x_n\}$ be a normalized unconditional basis of $X^*=Z_{1,p}$, $1<p<\infty$, (resp. $X^*=Z_{q,1}$, $1<q<\infty$) such a basis is boundedly complete and its biorthogonal basis spans a space isomorphic to $X=Z_{\infty,p^\prime}$ ( resp. $X=Z_{q^\prime,\infty}$).

We begin with a proposition which replaces Propositions \ref{pr:main1} and \ref{pr:main2} for the current cases.
\begin{pr}\label{pr:main4} Let $\{x_n\}$ be a normalized unconditional basis of $Z_{1,p}$, $1<p<\infty$, (resp. $Z_{q,1}$, $1<q<\infty$). Then $\{x_n\}$ contains a subsequence equivalent to the unit vector basis of $\ell_p$ (resp. $\ell_1$).
\end{pr}

\pf By proposition \ref{pr:main2}, $\ell_p$ does not embed into $Z_{1,p:\{n\},2}$ for $1<p<\infty$ and  $\ell_1$ does not embed into $Z_{q,1:\{n\},2}$ for $1<q<\infty$. It is thus enough to show that if $\{x_n\}$ contains no subsequence equivalent to the unit vector basis of $\ell_p$ (resp. $\ell_1$) then $[x_n]$ embeds in $Z_{1,p:\{n\},2}$ (resp. $Z_{q,1:\{n\},2}$).

The case of $Z_{q,1}$, $1<q<\infty$: We proceed as in the proof of Proposition \ref{pr:main1}. Since $q>1$ the beginning of the proof works for $p=1$ as well. The problem arise when we need to show that $\bar P$ is bounded as this no longer follow from Lemma \ref{lm:complementation}. But here we can use instead \cite[Theorem 1.d.6(ii)]{lt2} to prove that $\bar P$ is bounded in a very similar way to the proof of Lemma \ref{lm:complementation}. The rest of the proof of Proposition \ref{pr:main1} carries over.

The case of $Z_{1,p}$, $1<p<\infty$: Assume $\{x_n\}$ be a basis of $Z_{1,p}$, $1<p<\infty$. Let $\{x_n^*\}$ be the biorthogonal basis (of $Z_{\infty,p^\prime}$). By the assumption that $\{x_n\}$ doesn't
contain a subsequence equivalent to the unit vector basis of $\ell_p$, $[x_n^*]$ doesn't contain a subsequence equivalent to the unit vector basis of $\ell_{p^\prime}$.
The proof of Proposition \ref{pr:main1} works for $Z_{\infty,p^\prime}$, $1<p^\prime<\infty$, as well, with the same modification for the proof that $\bar P$ is bounded as in the previous paragraph, to show that in this case $[x_n^*]$ embeds (even complementably) into $Z_{\infty,p^\prime:\{n\},2}$.
\endpf

The next proposition replaces Proposition \ref{pr:main3} in the non-reflexive case.

\begin{pr}\label{pr:main5}(i) Let $\{x_n\}$ be a normalized unconditional basis of $Z_{1,p}$, $1<p<\infty$. Then the unit vector basis of $\ell_1$ is equivalent to a subsequence of $\{x_n\}$.

(ii) Let $\{x_n\}$ be a normalized unconditional basis of $Z_{q,1}$, $1<q<\infty$. Then the unit vector basis of $\ell_q$ is equivalent to a subsequence of $\{x_n\}$.
\end{pr}

\pf The proof of Proposition \ref{pr:main3} works for $Z_{q,p}$ also in the case $q=1<p<\infty$ and we get that under the assumption of $(i)$, if no subsequence of $\{x_n\}$ is equivalent to the unit vector basis of $\ell_1$ then $[x_n]$ embeds into $Z_{p,2}$ but this space has type $\min\{p,2\}$ so $\ell_1$ and thus also $Z_{1,p}$, $1<p<\infty$, do not embed into it. This proves $(i)$.

$(ii)$ It is enough to show that the unit vector basis of $\ell_{q^\prime}$ is equivalent to a subsequence of $\{x_n^*\}$ (the biorthogonal basis to $\{x_n\}$) which is an unconditional basis of $Z_{q^\prime,\infty}$. The proof of Proposition \ref{pr:main3} gives that if this is not the case then  $Z_{q^\prime,\infty}$ isomorphically embeds as a complemented subspace in $Z_{\infty,2}$.   Now if $Z_{q^\prime,\infty}$ isomorphically embeds as a complemented subspace in $Z_{\infty,2}$ then an easy application of Pe{\l}czynski's decomposition method gives that $Z_{q^\prime,\infty}\oplus Z_{\infty,2}$ is isomorphic to $Z_{\infty,2}$ but this immediately presents an unconditional basis for $Z_{\infty,2}$ which is not equivalent to a permutation of the cannonical basis of $Z_{\infty,2}$. This stands in contradiction to a result from \cite{bclt} and thus proves $(ii)$.
\endpf

%
%

\noindent G. Schechtman\\
Department of Mathematics\\
Weizmann Institute of Science\\
Rehovot, Israel\\
{\tt gideon@weizmann.ac.il}
\\

\end{document}